\documentclass[12pt]{article}

\usepackage{amssymb, amsthm, amsgen, amscd, amssymb, amsmath, amsxtra}

\usepackage{graphics, epsfig, color}

\usepackage[matrix, arrow, curve]{xy}

\usepackage{amsmath,amscd}

\newtheorem{defn}{Definition}

\newtheorem{thm}{Theorem}
\theoremstyle{remark}

\def\bee{\begin{equation*}}

\def\eee{\end{equation*}}

\def\be{\begin{equation}}

\def\ee{\end{equation}}

\def\C{{\mathbb C}}

\def\R{{\mathbb R}}

\def\Z{{\mathbb Z}}

\def\CP{{{\mathbb C}P}}

\title{The number of smooth $4$-manifolds with a fixed complexity}

\author{Dave Auckly
\thanks{This work was partially supported by National Science Foundation
grant DMS-0604994.}}

\begin{document}
\bibliographystyle{plain}
\maketitle

\centerline{Department of Mathematics, Kansas State University,
Manhattan, KS 66506 USA}

\centerline{\tt{dav@math.ksu.edu}}

\section{Introduction}
Donaldson demonstrated that topological $4$-manifolds could support
different smooth structures \cite{don, don2}. After this discovery
it was natural to ask how many different smooth structures a given
smooth $4$-manifold could support. For higher-dimensional manifolds
the answer was known -- a given (compact) topological manifold of
dimension $5$ or higher could only support a finite number of
different smooth structures. Quickly it was discovered that
topological $4$-manifolds could admit an infinite number of smooth
structures \cite{ok,fm}. Since then the techniques have improved and
we now know many examples of topological $4$-manifolds admitting
infinitely many different smooth structures \cite{FS}.

An alternate way to address the question is to ask how the number of
smooth structures grows with the complexity of the $4$-manifold.
Notions of complexity have been used in other settings in
low-dimensional topology \cite{mat}, and a notion of complexity for
$4$-manifolds was recently introduced in a nice paper of B. Martelli
\cite{mar}. Martelli proves that the number of homeomorphism classes
of simply-connected smooth $4$-manifolds with complexity less than
$n$ grows as $n^2$. This paper addresses the number of
diffeomorphism classes of simply-connected smooth $4$-manifolds with
complexity less than $n$, proving in particular that this number
grows at least as $n^{\sqrt[3]{n}}$. Along the way we construct
complete Kirby diagrams for a large family of knot surgery
manifolds.

Any smooth $4$-manifold admits a handle decomposition. A diagram
displaying the attaching regions is called a Kirby diagram; see
\cite{GS}. Roughly the complexity of a handle decomposition is the
sum of the number of disks, strands and crossings in the Kirby
diagram.  The complexity of a $4$-manifold is the minimal complexity
of a handle decomposition of the $4$-manifold. More precisely,
handlebodies are defined recursively with the empty set as the
simplest handlebody. A $k$-handle is a copy of $D^k\times D^{n-k}$
attached to a handlebody along the so-called attaching region
$S^{k-1}\times D^{n-k}$. The result of attaching a handle to a
handlebody is a new handlebody. It is not difficult to show that
every connected, smooth $4$ admits a handle decomposition with
exactly one $0$-handle. (On the other hand there are topological
$4$-manifolds that do not admit any handle decomposition.) The
attaching regions of the various handles can be put into general
position on the boundary of a unique $0$-handle. Since this boundary
is $S^3$ one can assume that the attaching regions of the
$1$-handles and $2$-handles miss one point. Removing this point
produces a copy of $\R^3$. Each $1$-handle will be attached along a
pair of $3$-disks. The cores of the $2$-handles intersect the
boundary of the $0$-handle in a compact $1$-manifold. (The core of
$D^k\times D^{n-k}$ is $D^k\times\{0\}$.) One can then take a
regular projection of this $1$-manifold to a plane. The result is a
Kirby diagram; see the figures in section \ref{diag} for examples.
The disks are the components of the attaching regions of the
$1$-handles. The strands are the components of the intersection of
the cores of the $2$-handles with the boundary of the $0$-handle and
the crossings are the crossings in the regular projection.

In order to obtain a lower bound on the growth of the number of
smooth structures one must first construct an interesting collection
of smooth $4$-manifolds and compute invariants to show that they are
distinct. Next one must construct Kirby diagrams for the manifolds,
and finally one will have to do a bit of combinatorics to estimate
the number as function of the complexity. This exactly outlines our
paper.

Fintushel and Stern gave a way to construct a $4$-manifold from a
knot, and related the Seiberg-Witten invariants of the $4$-manifold
to the Alexander polynomial of the knot \cite{FS}. Levine
constructed a family of knots producing every possible Alexander
polynomial \cite{le} completing the first part of the outline. These
constructions are reviewed in  section \ref{man} below. Akbulut and
Auckly independently described handle decompositions of the knot
surgery manifolds \cite{ak, au}. These decompositions are reviewed
and simplified in section \ref{diag} leaving a bit of combinatorics
for section \ref{cx} at the end of the paper.

I owe the referee thanks for very helpful comments that pointed out
errors in an earlier version of this paper.

\section{The manifolds}\label{man}
In order to build a family of smooth $4$-manifolds one should start
with one $4$-manifold. We start with an elliptic K3 surface. To be
precise let $\CP^3:=(\C^4-\{0\})/(\C-\{0\})$ be complex projective
space and define
\[
X:=\{[z_0:z_1:z_2:z_3]\in\CP^3|(z_0+z_1)^4-z_1^4+(z_2+z_3)^4-z_3^4=0\}\,.
\]
This is the Fermat surface. A simple application of the implicit
function theorem proves that it is a smooth $4$-manifold. Taken with
the projection $\pi:X\to \CP^1$ generically given by
$[z_0:z_1:z_2:z_3]\mapsto [z_0:z_2]$, this is an elliptic fibration.
The inverse image of $T:=\pi^{-1}([2:1])$ is
\[
\{[2z_2:z_1:z_2:z_3]\in\CP^3 |
17z_2^3+32z_2^2z_1+24z_2z_1^2+8z_1^3+4z_2^2z_3+6z_2z_3^2+4z_3^3=0\}\,.
\]
Further applications of the implicit function theorem demonstrate
that this fiber is a smooth, embedded, $2$-dimensional torus with
tubular neighborhood isomorphic to $T^2\times B^2$.

Using this torus we can apply the Fintushel-Stern knot surgery
construction to obtain a large family of homeomorphic $4$-manifolds.
Recall how this construction goes. Starting with a knot $K$ in the
$3$-sphere one identifies the boundary of $S^1\times
(S^3-\text{int}(N(K)))$ with $T^3$ so that $\text{pt}\times
S^1\times\text{pt}$ is a meridian of the knot and
$\text{pt}\times\text{pt}\times S^1$ is a longitude of the knot in
$\text{pt}\times S^3$. Finally one defines the knot surgery manifold
as
\[
X_K:=(X-T^2\times B^2)\cup_{T^3}(S^1\times (S^3-\text{int}(N(K)))\,.
\]
The powerful result from \cite{FS} is that the Seiberg-Witten
invariant of this manifold is given by
\[
\text{SW}_{X_K}=\Delta_K(t)\,.
\]
Here $\Delta_K$ is the Alexander polynomial of $K$ and
$t=\exp(2\text{PD}[T])$. The main point is that knots with different
Alexander polynomials give rise to distinct smooth $4$-manifolds.

\begin{figure}
\hskip50bp \epsfig{file=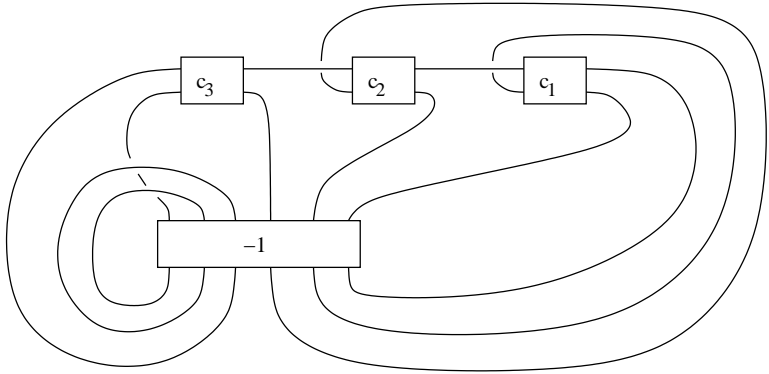, width=4truein} \caption{The knot
$K(c_1,c_2,c_3)$}\label{kc}
\end{figure}

To go further we need a family of knots with interesting Alexander
polynomials. We take the family displayed in figure \ref{kc}.
\begin{defn}
The Levine knots are the knots obtained by generalizing the knot
depicted in figure \ref{kc} to one with $d$ twist boxes along the
top, so that the number of strands passing down through the $(-1)$
twist box is the same as the number that pass up through the twist
$(-1)$ twist box. These knots will be denoted by
$K(c_1,c_2,\dots,c_d)$. We will denote the resulting knot surgery
manifold by $X(c_1,c_2,\dots,c_d):=X_{K(c_1,c_2,\dots,c_d)}$.
\end{defn}
This family was originally constructed by J. Levine to characterize
all possible Alexander polynomials \cite{le}.
Setting $c_0:=-1-2\sum_{k=1}^d c_k$ one can compute that
\[
\Delta_{K(c_1,c_2,\dots,c_d)}(t)=c_0+\sum_{k=1}^d c_k(t^k+t^{-k})\,.
\]
If the $(-1)$ twist box in the definition of $K(c_1,c_2,\dots,c_d)$
is changed to a $(+1)$ twist box then the value of $c_0$ would
change to $c_0:=1-2\sum_{k=1}^d c_k$ and the formula for the
Alexander polynomial would remain valid. This produces all possible
Alexander polynomials. It will be apparent from the Kirby diagrams
that the knot surgery manifolds $X_{-K}$ and $X_K$ are
diffeomorphic, so for the purposes of this paper there is no loss of
generality in using the $-1$ twist box. An exercise in Rolfsen's
book provides the computation of this Alexander polynomial
\cite{ro}. However, there is a sign error in the exercise so we give
a quick sketch of the computation.

\begin{figure}
\hskip50bp \epsfig{file=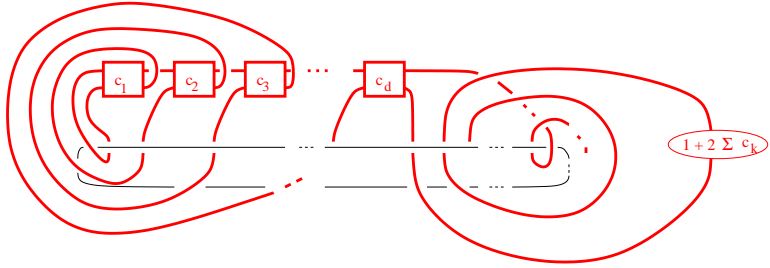, width=4truein} \caption{Surgery
description of the $K(c_1,\dots,c_d)$ complement}\label{isk}
\end{figure}

\begin{figure}
\hskip50bp \epsfig{file=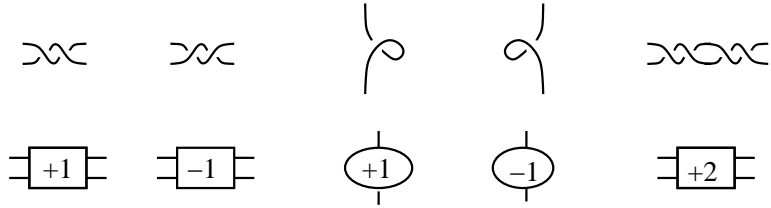, width=4truein} \caption{Twist
boxes and ovals}\label{twdef}
\end{figure}

First blow up the twist box by adding a $(+1)$-framed surgery
component. The resulting link can be isotoped into the projection
displayed in figure \ref{isk}.
In this figure the surgered component (one with twist boxes, red in
the electronic version) is drawn using the blackboard framing. To be
clear a $(+1)$-twist box, $(-1)$-twist box, $(+1)$-twist oval, and
$(-1)$-twist oval are displayed in figure \ref{twdef}. Higher order
twist boxes and ovals are defined as concatenations.

\begin{figure}
\hskip50bp \epsfig{file=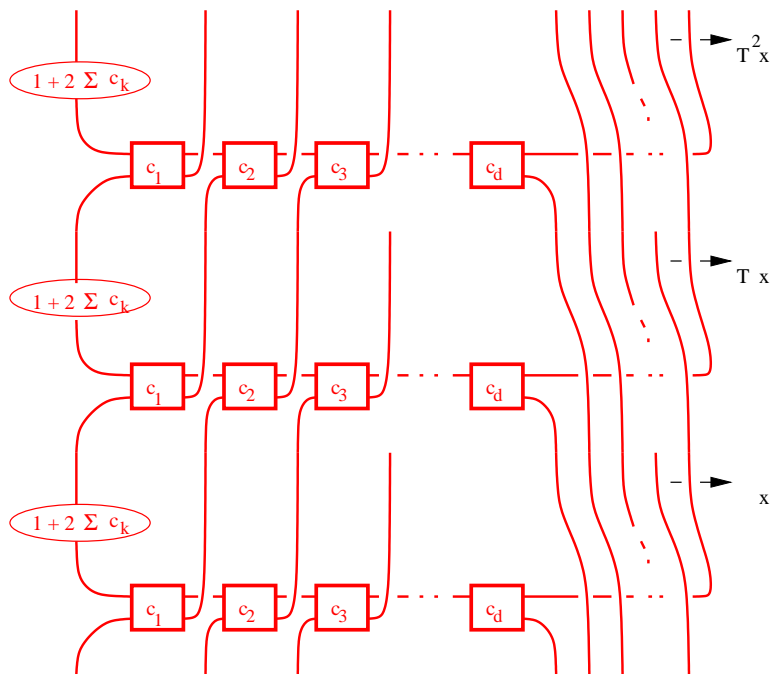, width=4truein} \caption{Infinite
cyclic cover of the $K(c_1,\dots,c_d)$ complement}\label{cov}
\end{figure}

It is easy to see that the infinite cyclic cover of the knot
complement has the surgery description displayed in figure
\ref{cov}. Clearly the first homology of the infinite cyclic cover
is generated by the $T^kx$ as an abelian group and is generated by
$x$ as a $\Z[T^{-1},T]$-module. The surgery curves will supply
relations. Following the surgery curve labeled with the $x$ from
just above the $x$ reading upward one can read off the relation.
Since there are $d-1$ crossings before a twist box is encountered,
the first twist box will be the $c_d$ box meeting the $T^dx$ surgery
curve. The next box will be the $c_{d-1}$ box on the $T^{d-1}x$
surgery curve, etc.. The last item encountered is the oval that
twists the framing. Each positive framing twist contributes a $-x$
to the relation for a total of $c_0=-1-2\sum_{k=1}^d c_k$. It
follows that the first homology of the infinite cyclic cover  is a
cyclic $\Z[T^{-1},T]$-module with relator
\[
\Delta_{K(c_1,c_2,\dots,c_d)}(t)=c_0+\sum_{k=1}^d c_k(t^k+t^{-k})\,.
\]

\nopagebreak
\section{Kirby diagrams}\label{diag}
In this section we construct a simple Kirby diagram for
$X(c_1,c_2,\dots,c_d)$. This starts with the procedure described in
\cite{ak,au}. Since there is a well-known Kirby diagram for the $K3$
surface (see \cite{GS}) we begin with a Kirby diagram for $S^1\times
(S^3-\text{int}(N(K(c_1,c_2,\dots,c_d))))$. Clearly this is a union
of two copies of $I\times
(S^3-\text{int}(N(K(c_1,c_2,\dots,c_d))))$. One obtains a properly
embedded $2$-disk in $D^4$ by taking the interval times the relative
pair obtained by removing a small ball containing a standard subarc
of any knot $K$. The boundary of the resulting $2$-disk is a copy of
$K\#-K$.
\begin{figure}
\hskip50bp \epsfig{file=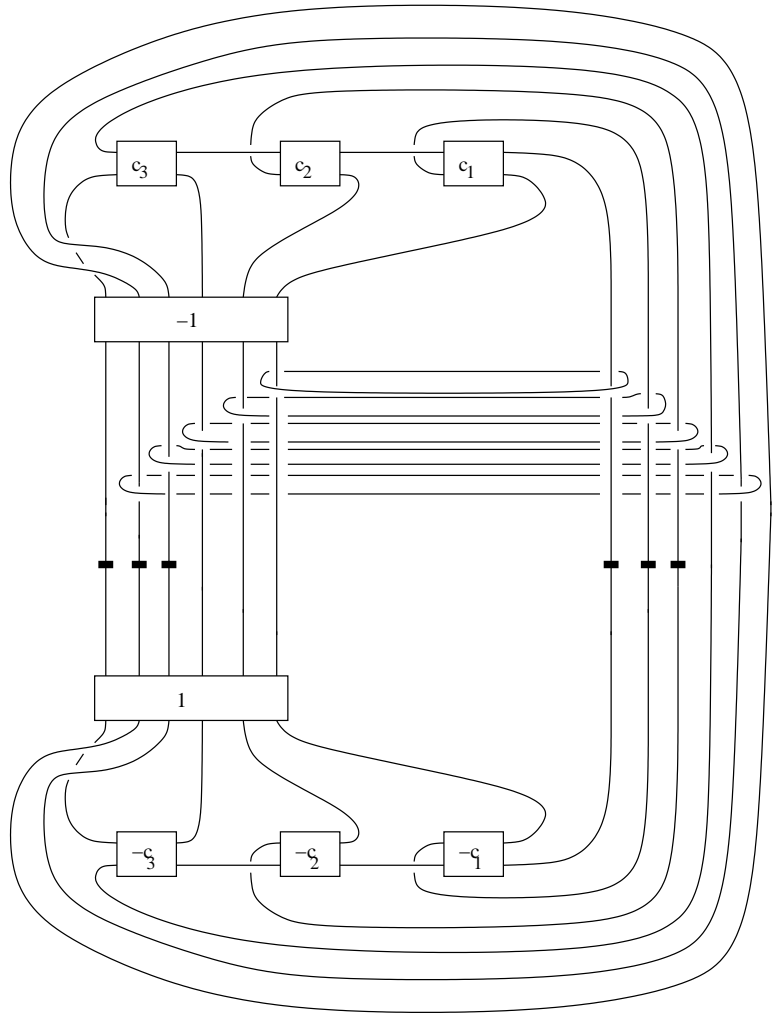, width=4truein} \caption{Complement
of the $K\#-K$ slice disk}\label{k-k}
\end{figure}
The complement of a tubular neighborhood of this $2$-disk is
described by putting a dot on the circle representing $K\#-K$. The
manifold obtained when this is applied to the unknot is easily seen
to be the same as the result of adding a $1$-handle to $D^4$. This
is a generalization of the standard `dotted circle' notation for
$1$-handles, and it generalizes in an obvious way to links.

In fact any properly embedded $2$-disk in $D^4$ in general position
with respect to the radius function produces a handle decomposition
of the complement of the disk with a $1$-handle in the complement
for each index $0$ critical point of the disk and a $2$-handle for
each index $1$ critical point of the disk etc.; see \cite{au}. The
result of this procedure applied to the disk obtained from
$K(c_1,c_2,\dots,c_d)$ is displayed in figure \ref{k-k}.

To double the handle decomposition here into one for \hfill\newline
$S^1\times (S^3-\text{int}(N(K(c_1,c_2,\dots,c_d))))$ notice that
doubling $I\times D^3$ amounts to adding a $1$-handle to $D^4$ and
to double the result of adding a $1$-handle to $I\times D^3$ amounts
to adding a $2$-handle etc. We call these the doubling handles. This
is explained in a bit more detail in \cite{au}. It follows that we
could complete the handle decomposition for $X(c_1,c_2,\dots,c_d)$
by adding one extra $1$-handle, the $2$-handles corresponding to the
$1$-handles in figure \ref{k-k}, the extra $2$-handles coming from
the decomposition of $K3$ and the $3$ and $4$-handles.

\begin{figure}
\hskip-2bp \epsfig{file=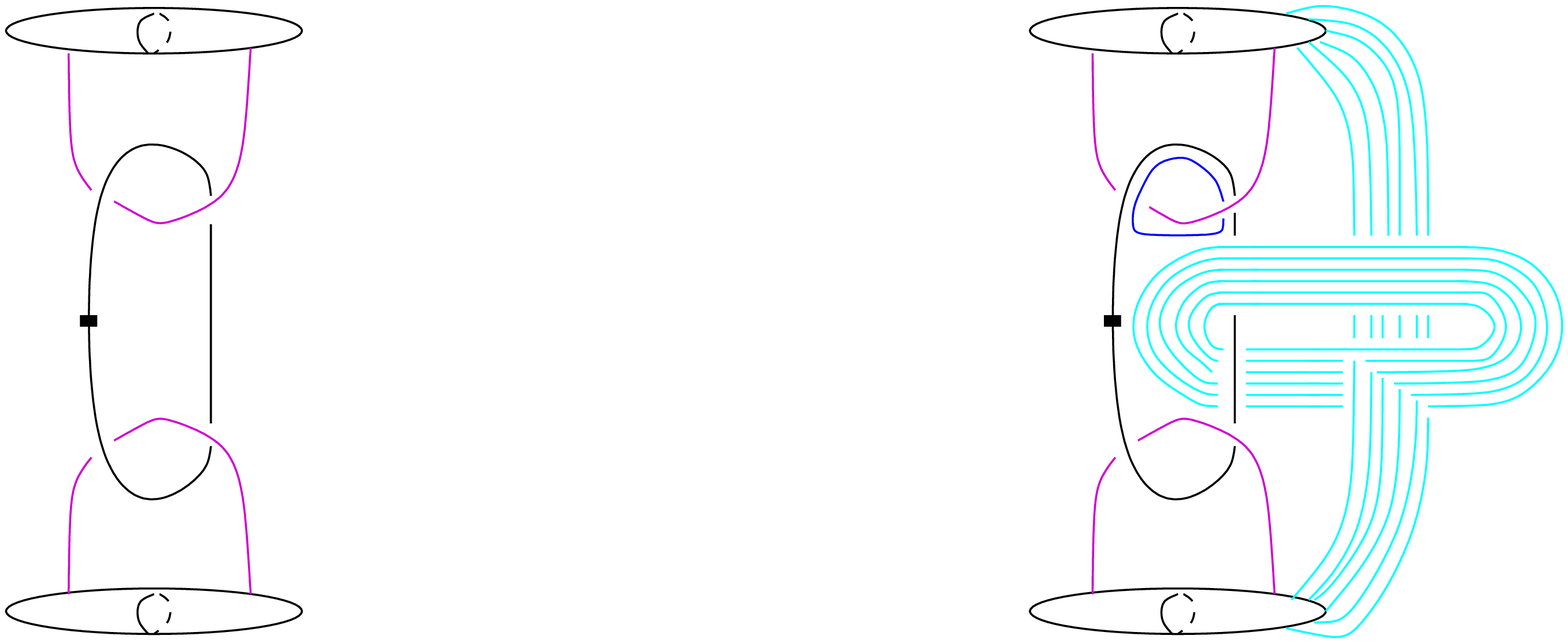, width=5.3truein} \caption{Torus
neighborhood and the rational elliptic surface}\label{k3}
\end{figure}

We first give the standard Kirby diagram for $K3$. Applying the
procedure that we described in this section to the unknot gives the
diagram for $T^2\times D^2$ on the left of figure \ref{k3}. The
$2$-handle in this diagram is taken with the $0$-framing.
The right side of figure \ref{k3} displays the Kirby diagram of the
rational elliptic surface. This figure may be obtained by drawing
the $1$-handle on the right of figure 8.17a of \cite{GS} as a footed
handle. The $2$ handles that were added to the $T^2\times D^2$ all
come with framing $-1$. The $2$-handle that is geometrically
unlinked from all of the $1$-handles is called the section handle.
The other twelve $2$-handles are from vanishing cycles. To obtain a
Kirby diagram for the $K3$ surface one just changes the framing on
the section to $-2$ and continues to add vanishing cycles in the
same pattern until there are a total of $24$ vanishing cycles.

\begin{figure}
\hskip50bp \epsfig{file=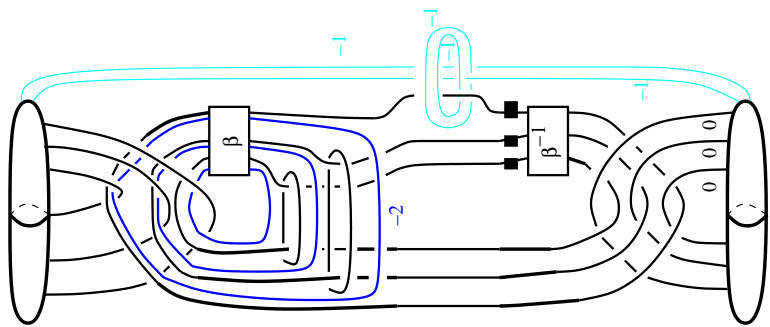, width=4truein} \caption{The
general knot surgery manifold}\label{knotsurg}
\end{figure}
To obtain a Kirby diagram for the knot surgery manifolds one just
needs to replace the $T^2\times D^2$ in the $K3$ diagram with the
diagram for $S^1$ times the knot exterior. See figure \ref{knotsurg}
for the general diagram obtained from this procedure. Here we are
assuming for simplicity that the knot is the closure of the braid
$\beta$. We further assume that the braid is drawn such that the
blackboard framing is the zero framing. This can always be done by
stabilizing the braid some number of times. Alternately the section
handle can be twisted around one of the $1$-handles until it
represents the longitude of the knot. The braid in the figure only
has three strands, but the generalization to higher order braids is
immediate. In addition this figure only includes four of the $24$
vanishing cycles.

The same procedure can be applied to any knot diagram, so we do not
have to worry that we don't have a braid presentation for our knot.
The result of this procedure applied to the knot $K(c_1,c_2,c_3)$
will have handles corresponding to doubling the thickened knot
complement, section, and vanishing cycles similar to figure
\ref{sch} but attached to the diagram of the thickened knot
complement from figure \ref{k-k}. The problem with this procedure in
our case is that we need to represent the $1$-handles by the feet of
the attaching regions. Thinking about the correspondence between the
dotted circle and $1$-handle one can see that strands linking the
dotted circle in the first representation correspond to strands that
pass over the handle in the second representation. It is also clear
that the collection of dotted circles representing the $1$-handles
must form an unlink. This is indeed the case for the diagram in
figure \ref{k-k}. However, isotoping the diagram until the dotted
circles are in standard position leads to many crossings and a
fairly complicated diagram. We will first do a bit of isotopy and
then add a few $1$/$2$-handle pairs in order to simplify the diagram
further. The result will be figure \ref{sch}.

From here forward our goal will be to generate a schematic of a
Kirby diagram with no dotted circles and fairly low complexity. A
complete Kirby diagram with no dotted circles would be very
complicated, but we only need to know enough to estimate the
complexity.

\begin{figure}
\hskip50bp \epsfig{file=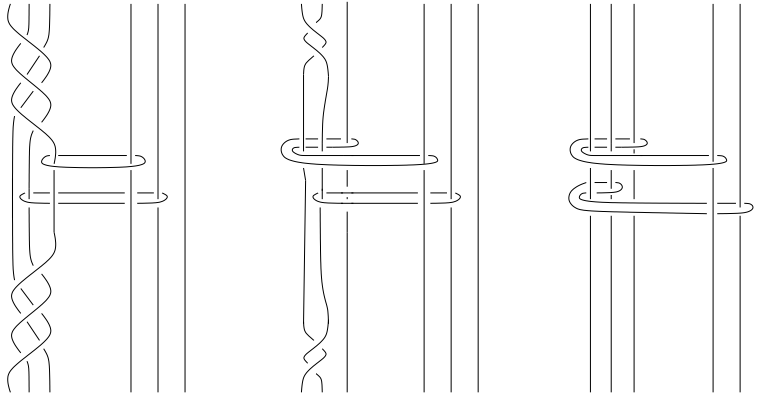, width=4truein}
\caption{Untwisting}\label{untwist}
\end{figure}

\begin{figure}
\hskip-2bp \epsfig{file=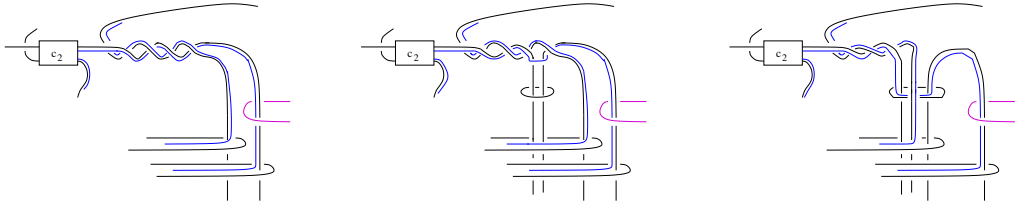, width=5.3truein} \caption{Adding
a $1$/$2$ pair}\label{pair1}
\end{figure}

The first isotopy will remove the pair of $(\pm 1)$-twist boxes as
in figure \ref{untwist}.
For the next part of the simplification we will keep track of some
of the $2$-handles from the $K3$ side as well. The elliptic
fibration $K3$ has a section. This is a $2$-sphere that meets each
fiber once. In particular when the $T^2\times B^2$ is removed a
$B^2$ is removed from the section, leaving a copy of $D^2$ in the
complement. Recall that the knot complement circle product is glued
so that a longitude of the knot is glued to the boundary of this
disk. A neighborhood of this disk becomes a $2$-handle attached to
the knot complement circle product. This is the $(-2)$-framed
$2$-handle in figure \ref{knotsurg}. This is why we call this
$2$-handle the section handle. Figure \ref{pair1} displays a portion
of the handle decomposition from figure \ref{k-k} after untwisting
together with part of the section handle and part of one of the
doubling $2$-handles, for $c_1=2$.

The pair of clasps on the bottom of the figure come from the topmost
$2$-handle from figure \ref{k-k} after untwisting. The clasp on the
right is part of the doubling $2$-handle associated to the
$1$-handle that it is linking. The loose strands that are mostly
parallel to existing $1$-handles are parts of the section handle. To
complete the figure one can add the reflection of all of the
$1$-handles  and the portion of the doubling $2$-handle across a
horizontal line at the bottom of the figure to the figure. This
would produce a larger part of the total Kirby diagram before adding
a handle pair. See figure \ref{sch} for the result after several
handle pair additions and isotopies.

The center diagram can be reflected in the same way without copying
the new $2$-handle to obtain the result of adding a $1$/$2$ pair. To
see that this is the result of adding a $1$/$2$ pair, slide the more
complicated $1$-handle passing through the left of the new
$2$-handle over the $1$-handle passing through the right and cancel
$1$/$2$ pair. The $2$-handle in the pair is just the meridian of the
dotted circle, so all handles that link the dotted circle can be
geometrically unlinked via this meridian. Notice that the section
handle may slide over the meridian $2$-handle, but its position does
not change in this process. Also, since the meridian is zero framed
and does not link any other $2$-handle the framings on all of the
$2$-handles are unchanged by this process.

\begin{figure}
\hskip-2bp \epsfig{file=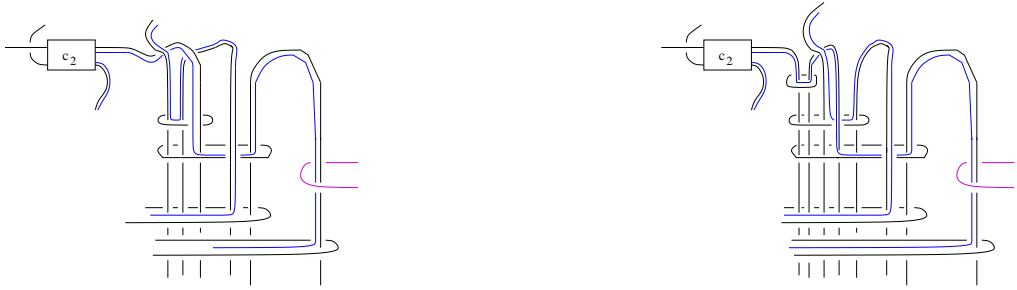, width=5.3truein}
\caption{Repeated handle additions}\label{pi}
\end{figure}

The diagram on the right shows the result after a bit of isotopy.
This procedure gets rid of a pair of crossings in the dotted circles
at a cost of adding a handle pair and a few more crossings. When we
started the section handle entered the crossing of the twist box
from the inside. After adding a $1$/$2$ pair and a bit of isotopy
the section handle enters the next crossing from the outside. The
same procedure can be done when the section handle enters the
outside of a crossing. Repeating this procedure, as in figure
\ref{pi}, we can get rid of all of the crossings from the first
twist box. In total we add $2c_k-1$ $1/2$ pairs to remove a $c_k$
twist box when $c_k>0$. Notice that for positive $c_k$ two crossings
can be isotoped away for free, but this does not happen for negative
$c_k$. For $c_k<0$ we can perform a small isotopy (Reidemeister II)
to make the twist box look exactly like the figure on the left of
figure \ref{pair1} with all of the crossings in the twist region
reversed. This means that we must add $2|c_k|+1$ pairs when $c_k<0$.
Even though this procedure adds more handles, it is still more
efficient than unwinding the $1$-handles in the original diagram.

\begin{figure}
\hskip-2bp \epsfig{file=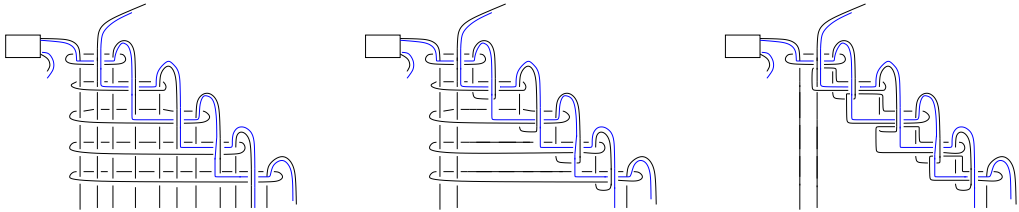, width=5.3truein} \caption{Twist box
after handle additions}\label{han}
\end{figure}

\begin{figure}
\hskip-2bp \epsfig{file=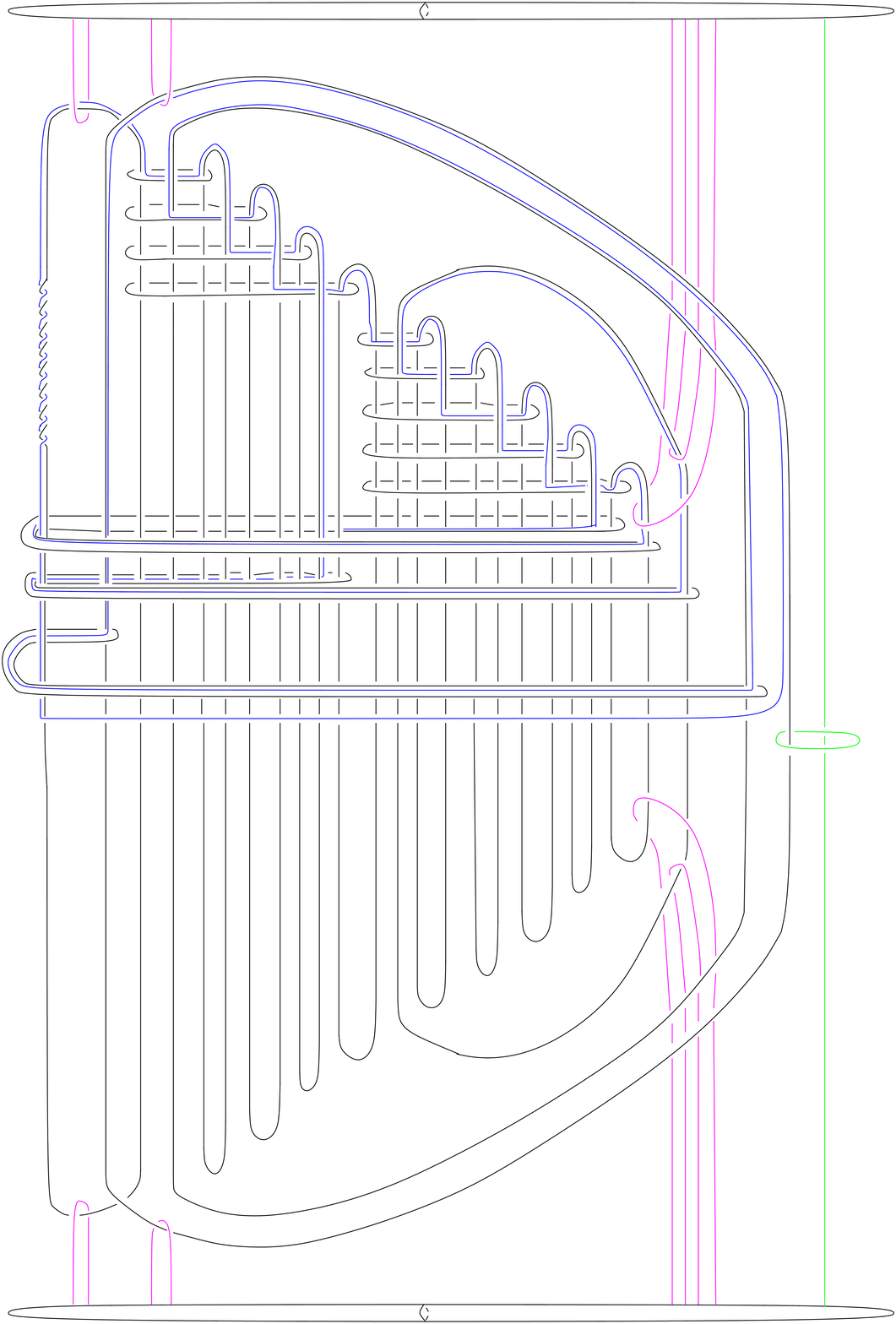, width=5.3truein} \caption{The
manifold $X(3,3)$ after some pair additions}\label{sch}
\end{figure}

Continuing to add handle pairs produces the diagram in figure
\ref{han}. (Here we draw the result obtained starting from $c_k=3$.)
The diagram in the center arises after a bit more isotopy. The
diagram on the right arises after a set of handle slides in which
the lowest $2$-handle is slid over the next highest until all of the
new $2$-handles take the `key' shape. The same procedure can be used
to remove the remaining twist boxes.
The result is displayed in figure \ref{sch} for $X(3,3)$, but this
figure only includes two of the twenty four vanishing cycles for
simplicity. Notice that the section handle forms a copy of the
original $K(3,3)$-knot.
The twists on the upper left arise because the section handle is
attached along a longitude to the original knot that does not link
the original knot. Referring to figure \ref{kc}, we see that each
twist box along the top contributes $-2c_k$ to the blackboard
framing, and the $(-1)$-twist box contributes $2d$ to the blackboard
framing. The remaining crossings come in pairs with opposite sign.
In order to have the blackboard framing equal to the zero framing,
we would have to add an oval representing a writhe of
$2\sum_{k=1}^dc_k-2d$. Also notice that we can isotope many of the
loops up and slide the new $2$-handles over each other as we did in
figure \ref{han}. In fact the number of loops that hang all the way
down from the remains of each twist box is independent of the number
of twists in the box.

\begin{figure}
\hskip-2bp \epsfig{file=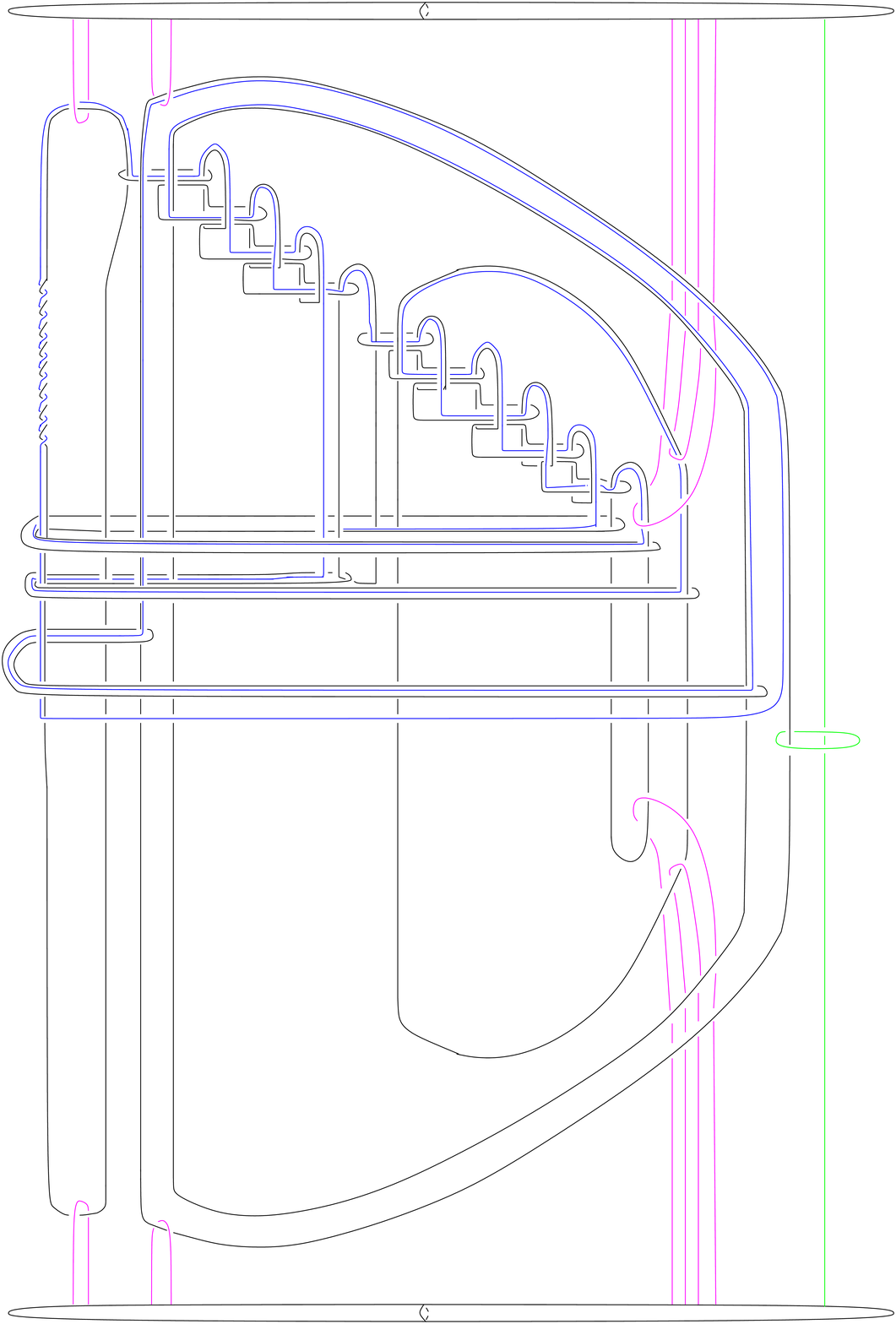, width=5.3truein} \caption{Final
Kirby diagram for $X(3,3)$}\label{sch2}
\end{figure}

The only remaining crossings between dotted circles are on the left
side of the diagram. These can be removed by a reasonable amount of
isotopy. In our diagram of $X(3,3)$ we first stretch the lowest hook
further, and then we can untwist. In general we have to stretch the
lowest $d-1$ hooks around in the same pattern and then we can
isotope.  The result is displayed in figure \ref{sch2}. The Kirby
diagram for $X(c_1,c_2,\dots,c_d)$ is just the natural
generalization of the diagram for $X(3,3)$.

\section{The complexity}\label{cx}

The Kirby diagram constructed in the last section is still not
suitable for a calculation of the complexity, but it is close. To
compute the complexity we need to draw the diagram without dotted
circles. We begin this section by describing how to convert a Kirby
diagram with knotted circles into one with $1$-handles.

\begin{figure}
\hskip32bp \epsfig{file=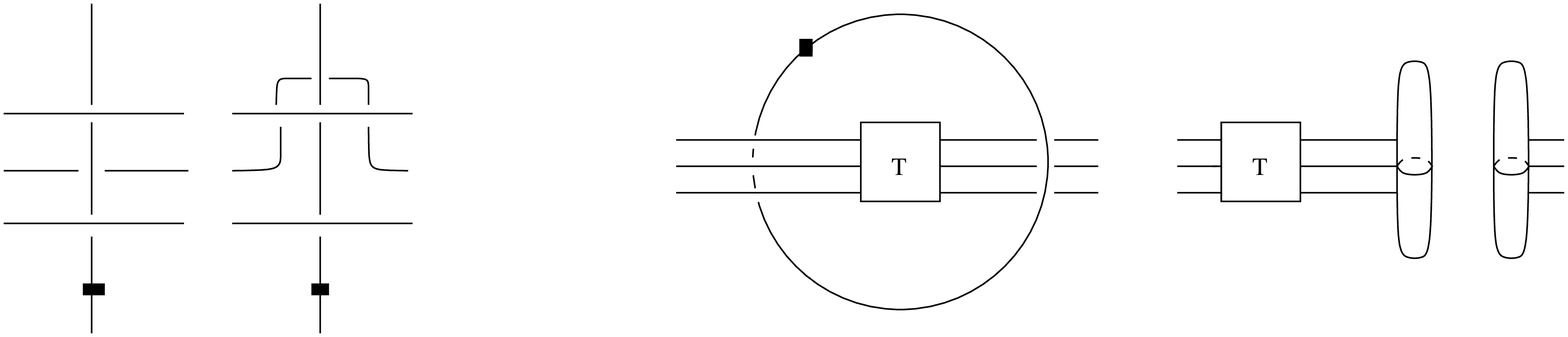, width=4.3truein} \caption{Changing
dotted circles into $1$-handles}\label{rm}
\end{figure}

The first step is to expand the dual handle decomposition of the
slice disks until the collection of dotted circles form an unlink
with no crossings. Each dotted circle will separate the plane of the
link projection into a bounded component and a non-bounded
component. The bounded component will contain a finite number of
circles and dotted circles together with a finite number of arcs
interacting in a tangle.  We would like all of the over-crossings to
be grouped together. Even though this might not be the case, we can
arrange for it to be the case by the move depicted on the left of
figure \ref{rm}. Once this is done the dotted circle in standard
position  can be converted into a $1$-handle as on the right of
figure \ref{rm}

What we need at this point is an upper bound on the complexity of
$X(c_1,c_2,\dots,c_d)$. Recall that the complexity of a manifold is
the minimal complexity of a Kirby diagram representing the manifold.
We will estimate the complexity of the diagram from figure
\ref{sch2}. Recall that the complexity is equal to the number of
disks plus the number of strands plus the number of crossings.
Clearly the number of disks is equal to twice the number of dotted
circles. The initial diagram for the knot complement times an
interval had $2d$ $1$-handles. We then added $2|c_k|+1$ $1$-handles
for each twist box (when $c_k$ is positive we could use $2$ fewer)
and there was one more from doubling. Putting this together shows
that
$$
\text{number of disks}\le 6d+4\sum_{k=1}^d |c_k|+2\,.
$$

Notice that the attaching circle of every $2$-handle meets at least
one $1$-handle. To count the number of strands it suffices to count
the number of under-crossings between $2$-handles and $1$-handles.
The $2$-handles coming from doubling the original $1$-handles
contribute $O(d^2)$ to this count. Since there is an upper-bound on
the number of crossings involving any given meridian $2$-handle from
the added handle pairs after the last isotopy, each contributes at
most a fixed number of strands and the number of these handles is
bounded by the total number of twists. Thus the $2$-handles in the
twist boxes contribute $O(\sum_{k=1}^d |c_k|+d)$ to the count (a
total of $d-1$ $1$-handles cross the furthest left such handle).
Each `hook' contributes $O(d)$ and there are $O(d)$ hooks for a
total of $O(d^2)$. The section handle contributes $O(\sum_{k=1}^d
|c_k|+d^2)$ to the count. This gives
$$
\text{number of strands}\le O(\sum_{k=1}^d |c_k|+d^2).$$

When we count the number of crossings we have to take into account
the fact that we have to do some moves to group all of the
under-crossings of any given dotted circle together, and these moves
produce extra crossings in the $2$-handles. There is an upper bound
on the number of crossings on each `small' $1$-handle in a twist
box, and there are $O(\sum_{k=1}^d |c_k|)$ of these `small'
$1$-handles. Ignoring the $1$-handle on the furthest left for the
moment, there are $O(d)$ remaining $1$-handles. Each of these
crosses $O(d)$ $2$-handles, but we might have to push as many as
$O(d)$ crossings past $O(d)$ strands as in the left of figure
\ref{rm}, so each of the remaining $1$-handles contributes $O(d^2)$
to the crossing count for a total of $O(d^3)$ crossings from these
handles. The $1$-handle on the far left has $O(|d-\sum_{k=1}^dc_k|)$
crossings coming from the writhe correction to the blackboard
framing. The over crossings from this area can be grouped together
using the move on the left of figure \ref{rm} resulting in a total
of $O((d-\sum_{k=1}^dc_k)^2)$ crossings. This handle also has $O(d)$
other crossings that can be grouped with the over and under from the
writhe area by the same moves contributing $O(d^2)$. Thus
$$ \text{number of crossings} \le O(d^3+\sum_{k=1}^d |c_k|)+O((d-\sum_{k=1}^dc_k)^2).$$ The
following estimate follows.
\begin{thm}
There are constants  $A_1$, $A_2$ and $A_3$ such that
$$
\text{\rm complexity}\left(X(c_1,c_2,\dots,c_d)\right)\le
A_1d^3+A_2\sum_{k=1}^d |c_k|+A_3\left(d-\sum_{k=1}^dc_k\right)^2\,.
$$
\end{thm}

We can turn this theorem around to obtain a lower bound on the
number of smooth classes of simply-connected $4$-manifolds with a
bounded complexity. First notice that
$$
\#\{(a_1,\dots,a_p)|a_k>0, \ \sum_{k=1}^pa_k=m\}=\binom{m-1}{p-1}\,,
$$
as there are $m-1$ `gaps' in a line of $m$ dots and such a sum
corresponds to choosing $p-1$ of the gaps. We can generalize this to
estimate the quantity below. By taking the first $[d/2]$ terms
positive and the rest to be negative and taking the sum of the
positive terms to be $\frac12(m+d)$ we obtain the lower bound
$$
\begin{aligned}
\#\{&(c_1,\dots,c_d)\in\Z^d|\sum_{k=1}^d|c_k|\le m, \ c_k\ne 0, \ \sum_{k=1}^dc_k=d \} \\
&\ge  \#\{(z_1,\dots,z_{[d/2]})|z_k>0, \ \sum_{k=1}^{[d/2]}
z_k=\frac12(m+d-2)\}\\ &\cdot \#\{(w_1,\dots,w_{d-[d/2]})|w_k>0, \
\sum_{k=1}^{d-[d/2]} w_k=\frac12(m-d-2)\} \\
&\ge \binom{\frac12(m+d-4)}{[d/2]}
\binom{\frac12(m-d-4)}{d-[d/2]}\,.
\end{aligned}
$$
Here this estimate is valid even if $m$ is not congruent to $d$ mod
$2$ or even not an integer (this is why we included the extra $-2$
terms).

Continuing, we see that
$$
\begin{aligned}
\#\{&(c_1,\dots,c_d)|
\text{complexity}\left(X(c_1,c_2,\dots,c_d)\right) \le n \} \\
&\ge \#\{(c_1,\dots,c_d)|
\text{complexity}\left(X(c_1,c_2,\dots,c_d)\right) \le n, \ \sum_{k=1}^dc_k=d \} \\
&\ge \#\{(c_1,\dots,c_d)| \sum_{k=1}^d|c_k|\le (n-A_1d^3)/A_2, \
c_k\ne 0, \ \sum_{k=1}^dc_k=d \}\,.
\end{aligned}
$$
Now pick $d=2[\frac14(n/A_1)^{1/3}]$ and continue the estimate for
large $n$

$$
\begin{aligned}
 &\ge
\binom{A_2^{-1}n/2-4A_1A_2^{-1}[\frac14(n/A_1)^{1/3}]^3+[\frac14(n/A_1)^{1/3}]-2}{[\frac14(n/A_1)^{1/3}]}
\\ &\cdot
\binom{A_2^{-1}n/2-4A_1A_2^{-1}[\frac14(n/A_1)^{1/3}]^3-[\frac14(n/A_1)^{1/3}]-2}{[\frac14(n/A_1)^{1/3}]}\\
&\ge \binom{A_2^{-1}n/3}{A_1^{-1/3}n^{1/3}/5}^2 \ge
n^{c\sqrt[3]{n}}\,.
\end{aligned}
$$
To obtain the last estimate we use $\binom{n}{k}\ge (n/k)^k$.  We
summarize this in the following theorem.
\begin{thm}
There is a constant $c$ so that for large $n$ the number of
diffeomorphism classes of simply-connected manifolds homeomorphic to
$K3$ having complexity less than or equal to $n$ is at least
$n^{c\sqrt[3]{n}}$.
\end{thm}
This should be compared with the following result of Martelli.
\begin{thm}[Martelli]
The number of homeomorphism classes of simply-connected
$4$-manifolds having complexity less than or equal to $n$ is between
$(1/4)n^2$ and $(5/16)n^2$.
\end{thm}
It should come as no surprise that the number of diffeomorphisms
grows much faster (faster than any polynomial) than the number of
homeomorphism classes (quadratic).

A more careful analysis of the crossings showing that the number of
crossings is $O(d^2+\sum_{k=1}^d |c_k|)$ might be possible. This
estimate would show that the number of diffeomorphism types grows at
least as $n^{c\sqrt{n}}$. It is harder to imagine improving the
bound much more than that with the same techniques. Other
constructions lead to similar combinatorics: using other elliptic
surfaces just changes the number of vanishing cycles by a constant;
distinguishing the result of several knot surgeries on distinct
fibers is easiest if each occurs in a $c$-neighborhood in which case
the complexity would appear to grow quadratically in the number of
knot surgeries; link surgeries also appear to have similar
combinatorics.

Martelli showed that the number of diffeomorphism classes of smooth
$4$-manifolds having complexity no greater than $n$ grows no faster
than $n^{Cn}$ for some constant $C$. It may be that this is the
right growth rate for simply-connected diffeomorphism types. To
prove this one would need to find a considerable simplification of
the Kirby diagrams presented here or find a different family of
manifolds with simple Kirby diagrams. Alternately one could look for
new $4$-manifold invariants. There are many more knots than the ones
that we have considered; however, the ones that we considered take
every possible value of the Alexander polynomial, so it is
impossible to get a Seiberg-Witten invariant other than the ones we
have here with a single knot surgery on a fiber of a $K3$.

It is interesting to ask the similar questions for $4$-manifolds
with additional structure. For example to address symplectic
$4$-manifolds it is natural to consider knot surgery with fibered
knots. This is similar in spirt to the work of Baldridge and Kirk
addressing how large a symplectic $4$-manifold with given
fundamental group must be \cite{BK}.

While it is clear that it is possible to answer a number of
questions about the complexity of $4$-manifolds, much less is known
about $4$-manifolds in general than is known about $3$-manifolds.
Thus questions of complexity are still more relevant in
$3$-dimensions where more subtle questions can be addressed
\cite{DT}.

\bibliography{bibcmplx}

\end{document}